\documentclass[oneside,11pt,reqno]{amsart}
\usepackage{stmaryrd}
\usepackage{bbm}
\usepackage{mathrsfs}
\usepackage{latexsym,amsxtra,xypic}
\usepackage[dvips]{graphicx}
\usepackage{dsfont}
\usepackage[all]{xy}
\usepackage{amscd,graphics}
\usepackage{amsmath,amsfonts,amsthm,amssymb}
\usepackage{latexsym,amsmath}
\usepackage{graphicx,psfrag}

\newtheorem{thm}{Theorem}[section]
\newtheorem{lem}[thm]{Lemma}
\newtheorem{cor}[thm]{Corollary}
\newtheorem{pro}[thm]{Proposition}
\newtheorem{ex}[thm]{Example}
\newtheorem{rmk}[thm]{Remark}
\newtheorem{defi}[thm]{Definition}
\newtheorem{fact}[thm]{Fact}
\newtheorem{nota}[thm]{Notation}
\newtheorem{ques}[thm]{Question}

\title
{The A-polynomial and holonomy perturbations}

\author{Jianfeng Lin}

\begin{document}

\maketitle 

\begin{abstract}
Dunfield-Garoufalidis and Boyer-Zhang proved that the A-polynomial of a non-trivial knot in $S^{3}$ is non-trivial. In this paper, we use holonomy perturbations to prove the non-triviality of the A-polynomial for a non-trivial, null-homotopic knot in an irreducible 3-manifold. Also, we give a strong constraint on the A-polynomial of a knot in the 3-sphere.
\end{abstract}

\section{Introduction}

In \cite{CCGLS}, an algebraic curve $D_{N}$ was associated to a compact 3-manifold $N$ with a single torus boundary. We briefly review the construction: consider the variety $\chi(\pi_{1}(N))$ of characters of $SL(2,\mathds{C})$ representations of $\pi_{1}(N)$. We can restrict a representation to the boundary group $\pi_{1}(\partial N)$, which gives us a map $r$: $\chi(\pi_{1}(N))\rightarrow \chi(\pi_{1}(T^{2}))$. Choose a basis $B=\{\mathfrak{M},\mathfrak{L}\}$ for the boundary group. Given two nonzero complex numbers $m,l$, we can define a $SL(2,\mathds{C})$ representation $\gamma_{(m,l)}$ of $\pi_{1}(T^{2})$:
\begin{equation}\label{diagnol}
\gamma_{(m,l)}(\mathfrak{M})=\left( \begin{array} {lr}
 m & 0 \\
 0& m^{-1}
\end{array}\right),
\gamma_{(m,l)}(\mathfrak{L})=\left( \begin{array} {lr}
 l & 0 \\
 0& l^{-1}
\end{array}\right)
\end{equation}

  So we have a natural map $h:\mathds{C}^{*}\times\mathds{C}^{*}\rightarrow\chi(\pi_{1}(T^{2}))$ defined by  $h(m,l):=[\gamma_{(m,l)}]\in\chi(\pi_{1}(T^{2}))$.

Now we consider $h^{-1}(\overline{\text{Im}(r)})$. The following lemma is proved in \cite{DG}.

\begin{lem}\label{dimension}
This algebraic set $h^{-1}(\overline{\text{Im}(r)})\subset \mathds{C}^{*}\times\mathds{C}^{*}$ has no 2-dimensional irreducible components.
\end{lem}

The algebraic curve $D_{N}$ is defined to be the union of the $1$-dimensional components of $h^{-1}(\overline{\text{Im}(r)})$. Each irreducible component of $D_{N}$ is defined by a two variable irreducible polynomial. We multiply all these polynomials to get the A-polynomial of $N$, which we denote by $A_{N,B}$. (It depends on the choice of the basis $B$.)

If $K$ is a null-homologous knot in a three manifold $Y$, we can take $N$ to be the knot complement. There is a natural basis $\widetilde{B}$ for the boundary group: the meridian and the longitude (see Fact \ref{meridian and longitude}). The A-polynomial $A_{K}$ of the knot $K$ is defined to be $A_{N,\widetilde{B}}$.

In \cite{DG} and \cite{Boyer-Zhang}, Dunfield-Garoufalidis and Boyer-Zhang independently proved that:
\begin{thm}
[Dunfield-Garoufalidis, Boyer-Zhang \cite{DG},\cite{Boyer-Zhang}]\label{non-trivial in sphere} For a non-trivial knot $K$ in $S^{3}$, we have $A_{K}\neq l-1$.
\end{thm}
\begin{rmk}\label{reducible}
The reducible representations of the knot group always give a factor $l-1$ of the A-polynomial. The A-polynomial of a trivial knot is exactly $l-1$. Therefore, we say a knot $K$ has non-trivial A-polynomial if $A_{K}\neq l-1$.
\end{rmk}

Both Dunfield-Garoufalidis's and Boyer-Zhang's proofs make use of the following theorem.
\begin{thm}
[Kronheimer-Mrowka \cite{KM}] For a non-trivial knot $K$ in $S^{3}$, let $Y_{r}$ ($r\in\mathds{Q}$) be the manifold obtained by doing $r$-surgery along $K$. If $|r|\leq2$, then $\pi_{1}(Y_{r})$ admits some $SU(2)$ representation with non-cyclic image.
\end{thm}

This theorem is proved using techniques from $SU(2)$ gauge theory, namely holonomy perturbations. Holonomy perturbations were used by Floer in \cite{FL} to prove the surgery exact triangle. In fact, using holonomy perturbations, we can prove Theorem \ref{non-trivial in sphere} directly. Moreover, we can prove something stronger.

The following are the main theorems of this paper.
\begin{thm}\label{main theorem0}
Suppose $K$ is a null-homologous, non-trivial knot in an orientable, closed, irreducible 3-manifold $Y$. If $A_{K}=l-1$, then there exists a homomorphism $\rho: \pi_{1}(Y )\rightarrow SU(2)$ such that $\rho([K])=-1\in SU(2)$.
\end{thm}

\begin{thm}\label{main theorem1}
For a null-homotopic, non-trivial knot $K$ in any orientable, closed, irreducible 3-manifold $Y$, we have $A_{K}\neq l-1$ .
\end{thm}

In light of these two theorems, we raise the following question.
\begin{ques}
Does any null-homologous knot in an irreducible 3-manifold have non-trivial A-polynomial?
\end{ques}

If we restrict to knots in $S^{3}$, Kronheimer and Mrowka's argument in \cite{KM} implies that the zero set of A-polynomial intersects the unit torus in some arcs. We give the precise statement:

\begin{thm}\label{3-sphere case}
 Let $K$ be a non-trivial knot in $S^{3}$. Then given any complex number $l$ with $|l|=1$, we can find another unit length complex number $m$ such that $A_{K}(l,m)=0$.
\end{thm}

\begin{rmk}
Theorem \ref{3-sphere case} gives strong constraints on the A-polynomial: for a generic polynomial, its zero set intersects the unit torus at isolated points. This is ruled out by Theorem \ref{3-sphere case}.
\end{rmk}

In section 2, we review some preliminaries and basic constructions related to holonomy perturbations and the non-vanishing theorem of critical points. In section 3, we prove the main theorems and give some examples.

\bigskip\noindent\textbf{Acknowledgement} The author wishes to thank Nathan Dunfield, Hans Boden and Yi Ni for valuable  discussions and comments. The author is especially grateful to Ciprian Manolescu for inspiring conversations and helpful suggestions in writing this paper.
\section{Preliminaries}

In this section, we review some of the constructions related to holonomy perturbations. Most of the details can be found in \cite{KM} and \cite{FLOER}.

\subsection{Holonomy perturbation}
Let $M$ be a closed, oriented 3-manifold with $b_{1}(M)>0$. Consider a rank 2 unitary bundle $E$ over $M$ with non-torsion $c_{1}(E)$. Let $\mathfrak{g}_{E}$ be the bundle whose sections are traceless, skew-hermitian endomorphisms of $E$. Let $\mathcal{A}$ be the affine space of $SO(3)$ connections of $\mathfrak{g}_{E}$. Let $\mathcal{G}$ be the group of gauge transformations on $E$ with determinant 1.

Fix a reference connection $A_{0}$ on $\mathfrak{g}_{E}$. For any connection $A$ on $\mathfrak{g}_{E}$, denoting $A-A_{0}$ by $\omega$, we have the Chern-Simons functional:
$$CS:\mathcal{A}\rightarrow \mathds{R}$$
$$CS(A)=\frac{1}{4}\int_{X_{0}} Tr(2\omega\wedge F_{A_{0}}+\omega\wedge d\omega+\frac{2}{3}\omega\wedge\omega\wedge\omega)$$

The critical points of the Chern-Simons functional are the flat connections.

Floer defined a class of perturbations of the Chern-Simons functional as follows. We take a function $\phi: SU(2)\rightarrow\mathds{R}$ which is invariant under conjugation. It is uniquely determined by the even, $2\pi$-periodic function $f(x):=$
$\phi\left(\begin{smallmatrix}
 e^{ix} & 0 \\
 0&e^{-ix}
\end{smallmatrix}\right)$.  Let $D$ be a compact 2-manifold with boundary. Consider an embedding of $D\times S^{1}$ in $M$ such that $\mathfrak{g}_{E}$ is trivial over it. Fix a trivialization of $\mathfrak{g}_{E}$ over $D\times S^{1}$ and take a 2-form $\mu$ which is supported in the interior of $D$ and has integral 1. Using the trivialization, we can lift $A$ to a connection $\bar{A} $ on the trivialized $SU(2)$ bundle $\widetilde{P}$ over $D\times S^{1}$. Consider:
$$\Phi : \mathcal{A}\rightarrow\mathds{R}$$
\begin{equation}\label{perturbation}
 \Phi(A)=\int_{p\in D}\phi(Hol_{\{p\}\times S^{1}}(\bar{A}))\mu(p)\end{equation}

In this paper, we specialize $D$ to be a disk $D^{2}$ or an annulus $H$. In the case of the disk, let $z_{0}\in \partial D^{2}$ be a base point. We define the curves $\alpha$ and $\beta$ to be $\partial D^{2}\times \{0\}$  and $\{z_{0}\}\times S^{1}$ respectively. In the case of the annulus, we denote the two components of $\partial H$ by $c_{1},c_{2}$ and choose base points $z_{i}\in c_{i}$. Then for $i=1,2$, we define $\alpha_{i},\beta_{i}\subset \partial (H\times S^{1})$ to be $c_{i}\times \{0\}$ and $\{z_{i}\}\times S^{1}$ respectively.

Now we consider the perturbed Chern-Simons functional: $CS+\Phi:\mathcal{A}\rightarrow \mathds{R}$.

The critical points can be completely described by the following lemmas. The first was proved by Floer. The second follows from similar arguments. For completeness, we sketch the proof for the second lemma. See \cite{FLOER} for details.
\begin{lem}[Floer \cite{FLOER}]\label{critical point}
  If $D\cong D^{2}$ and $A\in \mathcal{A}$ is a critical point of $CS+\Phi$, then:\\
1) $A$ is flat on $M-(D^{2}\times S^{1}) $\\
2) After choosing a new trivialization of $\widetilde{P}$, $\text{Hol}_{\alpha}(\bar{A})=\left(\begin{smallmatrix}
 e^{i\tau} & 0 \\
 0&e^{-i\tau}
\end{smallmatrix}\right)$, $\text{Hol}_{\beta}(\bar{A})=\left(\begin{smallmatrix}
 e^{i\nu} & 0 \\
 0&e^{-i\nu}
\end{smallmatrix}\right)$
and  $\tau= f'(\nu)+2\pi\mathds{Z}$\\
\end{lem}

\begin{lem}\label{critical point(2)}
If $D\cong H$ and $A\in \mathcal{A}$ is a critical point of $CS+\Phi$, then:\\
1) $A$ is flat on $M-(H\times S^{1}) $\\
2) The $SU(2)$ connection $\bar{A}$ is reducible and we can choose a new trivialization of $\widetilde{P}$ such that:

 $$\text{Hol}_{\alpha_{1}}(\bar{A})=\left(\begin{array} {lr}
 e^{i\tau_{1}} & 0 \\
 0&e^{-i\tau_{1}}
\end{array}\right), \text{Hol}_{\alpha_{2}}(\bar{A})=\left(\begin{array} {lr}
 e^{i\tau_{2}} & 0 \\
 0&e^{-i\tau_{2}}
\end{array}\right)$$ $$ \text{Hol}_{\beta_{1}}(\bar{A})=\text{Hol}_{\beta_{2}}(\bar{A})=\left(\begin{array} {lr}
 e^{i\nu} & 0 \\
 0&e^{-i\nu}
\end{array}\right)$$
Moreover, we have: \begin{equation}\label{boundary condition}\tau_{1}-\tau_{2}= f'(\nu)+2\pi\mathds{Z}\end{equation}
\end{lem}

\begin{proof}
Let $P$ be the trivialized $SU(2)$ bundle over $H\times S^{1}$ and $\mathfrak{g}_{E}$ be the vector bundle associated to the adjoint representation. For $p\in H$, denote $\{p\}\times S^{1}$ by $\beta_{p}$. For a point $a=(p,q)\in H\times S^{1}$, we consider the holonomy of $\bar{A}$ along $\beta_{p}$. Let us denote it by $T_{\bar{A}}(a)$. This is an automorphism of the fiber of $P$ over $a$. Thus $T_{\bar{A}}$ defines a section of $Aut(P)$. Since $\phi$ is conjugation invariant, it defines a function on $Aut(P)$. Let $\phi'(T_{\bar{A}})$ denotes the gradient of this function, evaluated at the section $T_{\bar{A}}$. One can check that
$\phi'(T_{\bar{A}})$ defines a section of $\mathfrak{g}_{E}$. Now $A$ is a critical point of $CS+\Phi$. It is proved in \cite{FLOER} that:

\begin{equation}\label{critical point equation}
F_{\bar{A}}=\phi'(T_{\bar{A}})\mu.
\end{equation}

Using this, we can deduce that the holonomy along $\beta_{p}$ does not depend on $p$ and  the covariant derivative $\nabla_{\bar{A}}(T_{\bar{A}})=0$ (see \cite{FLOER}). If $T_{\bar{A}}(a)=\pm1$ for any $a\in H\times S^{1}$, then  $\phi'(T_{\bar{A}})\equiv0$. By equation (\ref{critical point equation}), $\bar{A}$ is flat. Since $\pi_{1}(H\times S^{1})$ is abelian, we see that $\bar{A}$ is reducible. Suppose $T_{\bar{A}}(a)\neq\pm1$ for some $a\in H\times S^{1}$. Since $\nabla_{\bar{A}}(T_{\bar{A}})=0$, we have $T_{\bar{A}}(a)\neq\pm1$ for any $a\in H\times S^{1}$. Thus the existence of this section tells us that  $\bar{A}$ is also reducible in this case.

  Because $\bar{A}$ is reducible, we can choose a suitable trivialization of $\widetilde{P}$ such that $\bar{A}=\sigma\cdot\left(\begin{smallmatrix}
 i & 0 \\
 0& -i
\end{smallmatrix}\right)$ for some one  form $ \sigma\in \Omega^{1}(H\times S^{1})$. Since the holonomy along $\beta_{p}$ does not depend on $p$, we can assume that $ \text{Hol}_{\beta_{p}}(\bar{A})=\left(\begin{smallmatrix}
 e^{i\nu} & 0 \\
 0&e^{-i\nu}
\end{smallmatrix}\right)$ for any $p$. In particular, we have $\text{Hol}_{\beta_{i}}(\bar{A})=\left(\begin{smallmatrix}
 e^{i\nu} & 0 \\
 0&e^{-i\nu}
\end{smallmatrix}\right)$ for $i=1,2$. Then formula (\ref{critical point equation}) becomes:
\begin{equation}
d\sigma=f'(\nu)\mu.
\end{equation}
By Stokes' theorem, we have $\int_{\alpha_{1}}\sigma-\int_{\alpha_{2}}\sigma=\int_{H\times\{0\}}d\sigma=f'(\nu)$.
It is easy to see that $\tau_{i}=\int_{\alpha_{i}}\sigma$. Therefore, we proved the lemma.
\end{proof}

The following theorem was first proved in \cite{KM} and \cite{KM3}.
\begin{thm}
[Kronheimer-Mrowka \cite{KM},\cite{KM3}]Let $M,E,\mathfrak{g}_{E}$ be defined as before. If $M$ admits an oriented, smooth taut foliation and is not $S^{2}\times S^{1}$, then  for any holonomy perturbation $\Phi$, the perturbed Chern-Simons functional $CS+\Phi$ over $E$ has at least one critical point.
\end{thm}

This highly non-trivial theorem is a corollary of several theorems. We only sketch the proof here. For details, see \cite{KM} and \cite{KM3}. If $M$ admits a smooth, orientable taut foliation, then by the work by Eliashberg and Thurston in \cite{Eliashberg} and \cite{Eliashberg-Thurston}, $M$ can be embedded in an admissible symplectic 4-manifold $X$. \cite{KM} proved that $X$ can be chosen to satisfy several good conditions which imply, by Feehan and Leness's work (see \cite{KM3} and \cite{FL}), that it satisfies Witten's conjecture relating Seiberg-Witten invariants and Donaldson invariants. Because $X$ is symplectic, its Seiberg-Witten invariants are non-trivial by Taubes's non-vanishing Theorem in \cite{Taubes}. From Witten's conjecture, it follows that the Donaldson invariants $D^{v}_{X}$ are non-trivial for any line bundle $v$ on $X$. Then by a standard stretching argument, the perturbed Chern-Simons functional has at least one critical point.

The following two theorems are proved by Gabai.
\begin{thm}
[Gabai \cite{Gabai}] \label{taut foliation}Let $M$ be a closed, orientable, irreducible 3-manifold. If $b_{1}(M)>0$, then $M$ admits an orientable continuous taut foliation. Moreover, the taut foliation can be chosen to be smooth  if $H_{2}(M,\mathds{Z})$ is not generated by embedded tori or spheres.
\end{thm}

\begin{thm}[Gabai \cite{Gabai}] \label{taut foliation 2} Suppose $K\subset S^{3}$ is a non-trivial knot and $Y_{0}$ is obtained by doing $0$-surgery along $K$. Then $Y_{0}$ admits an orientable taut foliation and is not $S^{2}\times S^{1}$.
\end{thm}

\begin{rmk}
When the genus of $K\subset S^{3}$ is $1$, Gabai's taut foliation on $Y_{0}$ may not be smooth. But we can still embed $M$ in an admissible symplectic 4-manifold. See the proof of Theorem 6.1 in \cite{KMOS}.
\end{rmk}

Combining these theorems, we get:

\begin{cor}\label{exist critical point}
Suppose $M$ satisfies any one of the following two conditions:
\begin{itemize}
\item  $M$ is irreducible with $b_{1}(M)>0$ and $H_{2}(M,\mathds{Z})$ is not generated by embedded tori or spheres;

\item  $M$ is the 0-surgery manifold for a non-trivial knot in $S^{3}$.
\end{itemize}
Then for any holonomy perturbation on $M$, the perturbed Chern-Simons functional has at least one critical point.
\end{cor}

\subsection{The $SU(2)$-Pillowcase}

Now suppose $K\subset Y$ is a null-homologous knot in an orientable closed 3-manifold $Y$. Let $N(K)$ denote an open tubular neighborhood of $K$. Unless otherwise stated, all homology groups are considered with $\mathds{Z}$ coefficients. Using the Mayer-Vietoris sequence, it is easy to prove:

\begin{fact}\label{meridian and longitude}
The meridian $m\subset\partial(Y-N(K))$ is a primitive, nontorsion element in $H_{1}(Y-N(K))$. Also, we can find a unique longitude $l$ such that $[l]=0\in H_{1}(Y-N(K))$ and $([m],[l])$ forms a basis for $H_{1}(\partial(Y-N(K)))$.
\end{fact}

Let $Y_{r}(r\in \mathbb{Q})$ be the manifold obtained by doing $r$-surgery along $K$. It is easy to see that $ H_{1}(Y_{0})=H_{1}(Y-N(K))$.

Suppose that the attached solid torus is $N_{0}\subset Y_{0}$ with core $c$. Then $Y_{0}-N_{0}\cong Y-N(K)$. Now we consider a trivialized rank 2 unitary bundle $E$ over $Y-N(K)$ (using basic obstruction theory, we can prove that any $SU(2)$-bundle over $Y-N(K)$ is trivial since $\pi_{1}(SU(2))=\pi_{2}(SU(2))=0$). Let $\mathfrak{g}_{E}$ be defined as before. Let $\mathcal{G}$ be the group of gauge transformations with determinant $1$ on $E$.  Consider a flat connection $A$ on $\mathfrak{g}_{E}$. Using the trivialization, we can lift $A$ to a connection on the trivialized $SU(2)$ bundle $\tilde{P}$ over $Y-N(K)$. Denote this connection by $\tilde{A}$.

We have the following standard lemma:
\begin{lem}\label{holonomy}
By taking the holonomy of $\tilde{A}$, we get a 1-1 correspondence between the flat connections on $\mathfrak{g}_{E}$ modulo the action of $\mathcal{G}$ and the homomorphisms $\rho_{A}:\pi_{1}(Y-N(K))\rightarrow SU(2)$ modulo conjugation.
\end{lem}

Because $m$ and $l$ are commutative as elements in $\pi_{1}(Y-N(K))$, we can gauge transform $\tilde{A}$ so that $\text{Hol}_{m}(\tilde{A})$ and $\text{Hol}_{l}(\tilde{A})$ are both diagonal.

Thus we can assume that \begin{equation}\label{diagnol}
\text{Hol}_{m}(\tilde{A})=\left( \begin{array} {lr}
 e^{i\theta_{A}} & 0 \\
 0&e^{-i\theta_{A}}
\end{array}\right),
\text{ Hol}_{l}(\tilde{A})=\left( \begin{array} {lr}
 e^{i\eta_{A}} & 0 \\
 0&e^{-i\eta_{A}}
\end{array}\right).
\end{equation}

\begin{defi}\label{pillowcase} For a null-homologous knot $K\subset Y$, we define the subset of $(\mathds{R}/ 2\pi\mathds{Z})^{2}$:
 $$R_{K}:= \{\pm(\theta_{A},\eta_{A})|A \text{ is a flat connection on }  \mathfrak{g}_{E}|_{Y_{0}-N_{0}}\}$$
 \end{defi}

\begin{rmk}
The $SU(2)$-pillowcase refers to the space $\text{Hom}(\pi_{1}(T^{2}),SU(2))/\text{conjugation}$. Points in the torus  $(\mathds{R}/ 2\pi\mathds{Z})^{2}$ determine diagonal representations of $\pi_{1}(T^{2})$ to $ SU(2)$, although there is still a remaining conjugation action $(a,b)\mapsto (-a,-b)$. Therefore, $(\mathds{R}/ 2\pi\mathds{Z})^{2}$ is actually the double branch cover of the $SU(2)$-pillowcase. In this paper, we work with $\mathds{Z}/2\mathds{Z}$ equivariant arguments in the torus.
\end{rmk}

After checking the construction carefully and using Lemma \ref{holonomy}, it's easy to see that $R_{K}$ can be also defined as:
\begin{equation}\label{boundary}
\{(\theta,\eta)|\exists \rho:\pi_{1}(Y-K)\rightarrow SU(2) \text{ s.t. }\rho(m)=\left(\begin{array} {lr}
 e^{i\theta} & 0 \\
 0&e^{-i\theta}
\end{array}\right); \rho(l)=\left(\begin{array} {lr}
 e^{i\eta} & 0 \\
 0&e^{-i\eta}
\end{array}\right) \}
\end{equation}
\begin{rmk}\label{close subset}
Because $\pi_{1}(Y-K)$ is finitely generated and $SU(2)$ is compact, the $SU(2)$ representation space of  $\pi_{1}(Y-K)$ is compact. Thus $R_{K}$ is a closed subset of $(\mathds{R}/ 2\pi\mathds{Z})^{2}$.
\end{rmk}
\begin{rmk}\label{irreducible}
If $\rho$ is a reducible $SU(2)$ representation of $\pi_{1}(Y-K)$, then $\rho(l)=1$. Thus the points in $R_{K}$ off the line $\eta=2k\pi$ give irreducible representations.
\end{rmk}

The case $Y\cong S^{3}$ of the following lemma is proved as Lemma 13 in \cite{KM}. Although the general proof is essentially the same, we give it here for completeness.
\begin{lem}\label{translation}
$R_{K}$ is invariant under the translation by $(\pi,0)$.
\end{lem}

\begin{proof}
 By Fact \ref{meridian and longitude}, there exists $H_{1}(Y-N(K))\rightarrow \mathds{Z}_{2}$ mapping $m$ to $-1$ and mapping $l$ to $1$. Considering the composition of $\pi_{1}(Y-N(K))\rightarrow H_{1}(Y-N(K))\rightarrow \mathds{Z}_{2}\rightarrow SU(2)$, we get $\rho_{0}:\pi_{1}(Y-N(K))\rightarrow SU(2)$ with the image in the center. For $\rho: \pi_{1}(Y-N(K))\rightarrow SU(2)$, we multiply it with $\rho_{0}$ to get another representation $\rho'$. We have $\rho(m)=-\rho'(m)$ and $\rho(l)=\rho'(l)$. Then we use formula (\ref{boundary}) to prove the lemma.
\end{proof}

\section{Proof of main theorem}

\subsection{The general case}

Assume $K\subset Y$ is a null-homologous, non-trivial knot in an irreducible 3-manifold. If $K$ is contained in a 3-ball in $Y$, then $K$ has non-trivial A-polynomial by Theorem \ref{non-trivial in sphere}. Therefore, we just consider the case when $K$ is not contained in a 3-ball. This implies that $Y-N(K)$ is irreducible and its boundary is incompressible. We take two copies of $Y-N(K)$ and denote them by $Y_{j}-N(K_{j})$ ($j=1,2$). We glue them onto a thickened torus $T\times [-1,1]$ using the same gluing map to get a irreducible closed manifold $M$. $M$ is just the double of the knot complement. It's easy to see that $b_{1}(M)>0$ and the meridian $m\subset T\times \{0\}$ is a non-torsion element in $H_{1}(M)$.

\begin{rmk}
We identified the torus $T$ with $\partial(Y-N(K))$. Thus it makes sense to talk about the meridian $m\subset T$ and the longitude $l\subset T$ by the Fact \ref{meridian and longitude}.
\end{rmk}

\begin{lem}\label{nonvanishing thurston norm}
The Poincar\'{e} dual of $[m]\in H_{1}(M)$ (denoted by P.D.$[m]$) can not be represented by union of embedded tori or spheres. Therefore, the manifold $M$ satisfies the conditions of Corollary \ref{exist critical point}.
\end{lem}

\begin{proof}
Since $M$ is irreducible, embedded spheres are all null-homologous. Thus we just consider the case $P.D.[m]=\sum [T_{i}]$ where $T_{i}$ are embedded tori. We can assume each $T_{i}$ intersect both $m$ and $T\times \{0\}$ transversely.  Since the intersection number of $P.D.[m]$ with $m$ equals 1, some $T_{i}$ must intersect $m$ for odd times. Without loss of generality, we can assume this torus to be $T_{1}$. Now suppose $T_{1}\cap( T\times \{0\})=\mathop{\cup}\limits_{j=1}^{n} \gamma_{j}$ where $\gamma_{j}$ are disjoint simple closed curves. We can assume none of these curves bound disks in $T_{1}$. Otherwise, by choosing the innermost curve which bounds a disk, we are able to find some $\gamma_{i}$ bounding a disk $D\subset T_{1}$ and the interior $\mathring{D}$ does not intersect $T\times \{0\}$. Since $T\times \{0\}$ is incompressible, $\gamma_{i}$ also bounds a disk $D'\subset T\times \{0\}$. Notice that $D\cup D'$ is null-homologous since $M$ is irreducible. We can replace $T_{1}$ with $(T_{1}\setminus D)\cup D'$ to eliminate the intersection $\gamma_{i}$.

Moreover, suppose some $\gamma_{i}$ bounds a disk $D''\in T\times \{0\}$. As before, we can assume that $\mathring{D}''$ does not intersect $T_{1}$ by choosing the innermost curve. Consider the surface $T_{1}'=((M\setminus (T\times[-1,1]))\cap T_{1})\cup (D''\times\{\pm1\})$. Notice that $\gamma_{i}$ is a non-separating curve in $T_{1}$ since it does not bound a disk. We see that $T_{1}'$ is an embedded sphere representing the same  (non-zero) homology class as $T_{1}$, which is a contradiction with the fact that $M$ is irreducible.

Now each $\gamma_{i}$ is not null-homologous in either $T_{1}$ or $T\times \{0\}$. Since these curves don't intersect each other, they are parallel in both $T_{1}$ and $T\times\{0\}$. Notice that the components of $T_{1}\setminus (\mathop{\cup}\limits_{j=1}^{n}\gamma_{j})$ are contained in the two components of $M\setminus (T\times \{0\})$ alternatively, which implies that $n$ is even. Since each $\gamma_{i}$ intersects $m$ for the same times, the intersection number of $T_{1}$ with $m$ is even. This is a contradiction.\end{proof}

Now we consider the rank 2 unitary $E$ with $c_{1}(E)$ the Poincar\'{e} dual of $[m]\in H_{1}(M)$. We can get $\mathfrak{g}_{E}$ by the following procedures:

\begin{itemize}
\item Construct trivialized $SO(3)$ bundles over $Y_{j}-N(K_{j})$ ($j=1,2$) and $T\times [-1,1]$. Denote them by $\mathfrak{g}_{1}$, $\mathfrak{g}_{2}$ and $\mathfrak{g}_{0}$ respectively.

\item Choose a map $f:T^{2}\rightarrow SO(3)$ such that  $f_{*}:\pi_{1}(\partial{N_{0}})\rightarrow \pi_{1}(SO(3))\cong\mathds{Z}_{2}$ satisfies:
 \begin{equation}\label{diff} f_{*}([m])=[0],  f_{*}([l])=[1]. \end{equation}
\item Glue $\mathfrak{g}_{1}$, $\mathfrak{g}_{0}$ along $T\times\{-1\}$ using the identity map and glue $\mathfrak{g}_{0}$,$\mathfrak{g}_{2}$ along $T\times \{+1\}$ using $f$.
\end{itemize}

In order to do holonomy perturbations, we embed $H_{j}\times S^{1} (j=1,2)$ into $M$ by $i_{1},i_{2}$ as follows:

$i_{1}(H_{1}\times S^{1})= T\times [-1,0]$ and $i_{1}(*\times S^{1})$ is parallel to $l\subset T^{2}\times [-1,0]$, which means that we use the holonomy along the longitude to do perturbation $\Phi_{1}$.

$i_{2}(H_{2}\times S^{1})= T\times [0,1]$ and $i_{2}(*\times S^{1})$ is parallel to $m\subset T^{2}\times [0,1]$, which means that we use the holonomy along the meridian to do perturbation $\Phi_{2}$.

Now we consider the perturbed Chern-Simons functional, $CS+\Phi_{1}+\Phi_{2}:\mathcal{A}\rightarrow \mathds{R}$.

Notice that we already have trivializations of $\mathfrak{g}_{E}$ when restricted to $Y_{j}-N(K_{j})$ and $T^{2}\times [-1,1] $. Thus for a connection $A\in \mathcal{A}$, we can restrict $A$ to these pieces and use the trivializations to lift it to trivialized $SU(2)$ bundles. We have six loops: $m\times \{j\}\subset T^{2}\times[-1,1]$ and $l\times \{j\}\subset T^{2}\times[-1,1], j=-1,0,1$. Denote them by $m_{j}, l_{j}$.

If $A$ is a critical point of $CS+\Phi_{1}+\Phi_{2}$, denote by $\widetilde{A}$ the $SU(2)$ lift of $A|_{\mathfrak{g}_{0}}$.

Since we only do perturbations in the interior of $T\times [0,-1]$ and $T\times [0,1]$, we have $\widetilde{A}|_{T\times \{j\}}$ flat for $j=-1,0,1$. Taking the holonomy of $\widetilde{A}|_{T\times \{j\}}$, we get a representation $\rho_{j}:\pi_{1}(T^{2})\rightarrow SU(2)$.
\begin{nota}
For $(\theta,\eta)\in (\mathds{R}/2\pi\mathds{Z})^{2}$, we denote by $\rho_{(\theta,\eta)}$ the representation $\pi_{1}(T^{2})\rightarrow SU(2)$ which maps $m$ to $\left( \begin{smallmatrix}
 e^{i\theta} & 0 \\
 0&e^{-i\theta}
\end{smallmatrix}\right)$ and  $l$ to  $\left( \begin{smallmatrix}
 e^{i\eta} & 0 \\
 0&e^{-i\eta}
\end{smallmatrix}\right)$. We write $\rho\sim\rho'$ if two representations $\rho$ and $\rho'$ differ by a conjugation. It is easy to see that $\rho_{(\theta,\eta)}\sim \rho_{(\theta',\eta')}$ if and only if $(\theta,\eta)=\pm(\theta',\eta')$.
\end{nota}

 Applying Lemma \ref{critical point(2)} to $\Phi_{1}$ with $\alpha_{1}=m_{-1}, \alpha_{2}=m_{0},\beta_{1}=l_{-1},\beta_{2}=l_{0}$, we see that $\rho_{-1}\sim \rho_{(\theta_{-1},\eta_{-1})}$ and
$\rho_{0}\sim \rho_{(\theta_{0},\eta_{0})}$ with:

 \begin{equation}\label{diff2}\eta_{-1}=\eta_{0} \text{ and } \theta_{-1}-\theta_{0}=f_{1}'(\eta_{0}).\end{equation}

 Applying Lemma \ref{critical point(2)} to $\Phi_{2}$ with $\alpha_{1}=l_{0}, \alpha_{2}=l_{1},\beta_{1}=m_{0},\beta_{2}=m_{1}$, we see that $\rho_{0}\sim \rho_{(\hat{\theta}_{0},\hat{\eta}_{0})}$ and
$\rho_{1}\sim \rho_{(\theta_{1},\eta_{1})}$ with:
 \begin{equation}\label{diff3}\hat{\theta}_{0}=\theta_{1} \text{ and } \hat{\eta}_{0}-\eta_{1}=f_{2}'(\theta_{1}).\end{equation}

\begin{rmk}
Recall that we choose $\phi_{1}, \phi_{2}: SU(2)\rightarrow \mathds{R}$ to define $\Phi_{1}$ and $\Phi_{2}$. They give us $f_{1}$, $f_{2}$.
\end{rmk}

We can also use the trivialization of $\mathfrak{g}_{j} (j=1,2)$ to lift $A$ to flat $SU(2)$ connections over $Y_{j}-N(K_{j})$. Denote them by $\widetilde{A}_{1}$ and $\widetilde{A}_{2}$. Because we glue $\mathfrak{g}_{1}$, $\mathfrak{g}_{0}$ using identity map, we have:$$\text {Hol}_{m_{-1}}(\widetilde{A})=\text {Hol}_{m_{-1}}(\widetilde{A}_{1}),\text {Hol}_{l_{-1}}(\widetilde{A})=\text {Hol}_{l_{-1}}(\widetilde{A}_{1})$$

 Since we glue $\mathfrak{g}_{0}$ and $\mathfrak{g}_{2}$ using $f$, their trivializations do not agree on their common boundary. Different trivializations of the $SO(3)$ bundle give different $SU(2)$ lifts of $A|_{T\times\{2\}}$. By formula (\ref{diff}), it is easy to check:

$$\text {Hol}_{m_{1}}(\widetilde{A})=\text {Hol}_{m_{1}}(\widetilde{A}_{2}),\text {Hol}_{l_{1}}(\widetilde{A})=-\text {Hol}_{l_{1}}(\widetilde{A}＆_{2})$$

\begin{nota}
Let $S\subset (\mathds{R}/2\mathds{Z}\pi)^{2}$ be a subset. If $h$ is a function with period $2\pi$, we denote the set $\{(\theta,\eta+h(\theta))| (\theta,\eta)\in S \}$ by $S+(*,h)$ and the set $\{(\theta+h(\eta),\eta)| (\theta,\eta)\in S \}$ by $S+(h,*)$.
We also denote the set $\{(\theta+a,\eta+b)| (\theta,\eta)\in S \}$ by $S+(a,b)$ for constant $a,b$.

\end{nota}

By Definition \ref{pillowcase} and the discussion above, we see that
\begin{equation}\label{critical point(4)}
(\theta_{-1},\eta_{-1})\in R_{K}, (\theta_{1},\eta_{1})\in R_{K}+(0,\pi).\end{equation}

By formula (\ref{diff2}), we have $(\theta_{0},\eta_{0})\in R_{K}+(-f_{1}',*)$. Since $\rho_{(\theta_{0},\eta_{0})}\sim \rho_{0}\sim \rho_{(\hat{\theta}_{0},\hat{\eta}_{0})}$, we have $(\hat{\theta}_{0},\hat{\eta}_{0})=\pm (\theta_{0},\eta_{0})$. Notice that $R_{K}+(-f_{1}',*)$ is symmetric under the reflection at $(0,0)$. We have $(\hat{\theta}_{0},\hat{\eta}_{0})\in R_{K}+(-f_{1}',*)$. By formula (\ref{diff3}), we see that $(\theta_{1},\eta_{1})\in (R_{K}+(-f_{1}',*))+(*,-f_{2}')$.  Thus $(\theta_{1},\eta_{1})$ is an intersection point of  $(R_{K}+(-f_{1}',*))+(*,-f_{2}')$ and $R_{K}+(0,\pi)$. We get the following important proposition.

\begin{pro}\label{non-existance}
Suppose that $(R_{K}+(-f_{1}',*))+(*,-f_{2}')$ does not intersect $R_{K}+(0,\pi)$. Then the perturbed Chern-Simons functional $CS+\Phi_{1}+\Phi_{2}$ has no critical point.
\end{pro}

Now we are ready to prove Theorem \ref{main theorem0}.
 \begin{figure}[htbp]
\centering
\includegraphics[angle=0,width=0.5\textwidth]{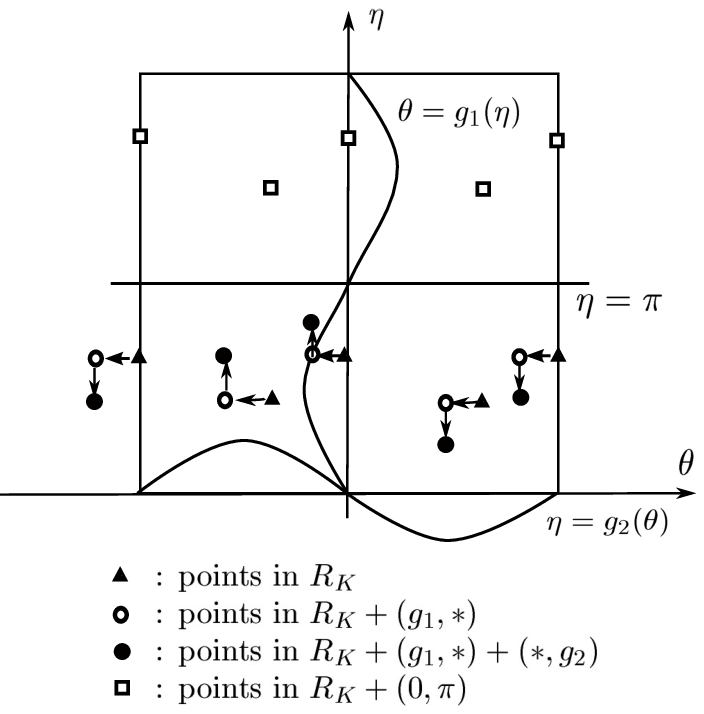}
\caption{} \label{meinv}
\end{figure}

 \begin{proof}

 Suppose $K$ has trivial A-polynomial. Then $R_{K}\setminus\{\eta=2k\pi\}$ has only finitely many points, by Lemma \ref{dimension}. If there exists no representation $\pi_{1}(Y)\rightarrow SU(2) $ mapping $K$ to $-1$, then $(0,\pi)\notin R_{K}$. By Lemma \ref{translation} $(k\pi,\pm\pi)\notin R_{K}$. So, we can choose an odd, $2\pi$-periodic function $g_{1}$ such that the curve $\{g_{1}(\eta)=\theta, 0<\eta<2\pi\}$ does not intersect $R_{K}$. Then $R_{K}+(g_{1},*)$ does not intersect $\{\theta=k\pi,0<\eta<2\pi\}$. This can be done because there are only finitely many points to move away from the lines $\{\theta=k\pi,0<\eta<2\pi\}$ and none of the them have $\eta$ component equal to $\pi$. We find another odd, $2\pi$-periodic  function $g_{2}$ such that $(R_{K}+(g_{1},*))+(*,g_{2})$ does not intersect the translation $R_{K}+(0,\pi)$. This can be done because $R_{K}+(g_{1},*)$ has only finitely many points outside the circle $\theta=2k\pi$ and none of them have $\theta$ component equal to $k\pi$. $R_{K}+(0,\pi)$ has only finitely many points off the circle $\theta=\pi$. We can move these discrete points away from the one dimensional components. As long as $|g_{2}|<\frac{\pi}{2}$, the one dimensional component of  $(R_{K}+(g_{1},*))+(*,g_{2})$ and that of $R_{K}+(0,\pi)$ do not intersect.

Then we find even, $2\pi$-periodic  functions $f_{1},f _{2}$ such that $f_{j}'=-g_{j}$ and use them to define perturbations $\Phi_{j}$. By Proposition \ref{non-existance}, the perturbed Chern-Simons functional $CS+\Phi_{1}+\Phi_{2}$ has no critical point.
However, since $K$ is not contained in a $3$-ball by Theorem \ref{non-trivial in sphere}, $CS+\Phi_{1}+\Phi_{2}$ has at least one critical point by Lemma \ref{nonvanishing thurston norm} and Corollary \ref{exist critical point}. This is a contradiction.
\end{proof}

Theorem \ref{main theorem1} is a trivial corollary of Theorem \ref{main theorem0}. Here are some other corollaries.

\begin{rmk}
If $K$ is null-homotopic, then it is easy to see that $R_{K}$ has no points on the line $\theta=0,0<\eta<2\pi$. Therefore, Theorem \ref{main theorem1} can be proved using the perturbation $\Phi_{2}$ without the use of $\Phi_{1}$.
\end{rmk}

\begin{cor}
If $K$ is a non-trivial, null-homologous knot in an irreducible, orientable 3-manifold $M$ and $\pi_{1}(M)/(K)$ has no $SO(3)$ representation with non-cyclic image, then $A_{K}\neq l-1$. Here $(K)$ denotes the normal subgroup group generated by $[K]\in \pi_{1}(M)$. In particular, if $\pi_{1}(M)$ has no non-cyclic $SO(3)$ representation, then any null-homologous, non-trivial knot has non-trivial A-polynomial.
\end{cor}

\begin{proof}
Suppose $K$ has trivial A-polynomial. Then there exists $\rho:\pi_{1}(M)\rightarrow SU(2) $ such that $\rho(K)=-1$. The representation $\rho$ induces a $SO(3)$-representation of $\pi_{1}(M)$ which factors through $\pi_{1}(M)/(K)$. By our assumption, this representation has cyclic image and so does $\rho$. Thus $\rho$ factors through $H_{1}(M)$, which is impossible because $K$ is null homologous and $\rho(K)\neq 1$.
\end{proof}

\begin{ex}
If $K\subset M$ is a non-trivial knot in the Poincar\'{e} homology sphere and $2[K]\neq e\in \pi_{1}(M)$, then $A_{K}\neq l-1$.
\end{ex}

\begin{ex}
Let $K\subset Y_{1}$ be a non-trivial knot in $Y_{1}$, where $Y_{1}$ is the manifold obtained by doing 1-surgery along some knot $K'$ in $S^{3}$. If $K$ is homotopic in $Y_{1}$ to the meridian of $S^{3}-K'$ then $A_{K}\neq l-1$.

\end{ex}

The following example is given by Dunfield:
\begin{ex}
Let $Y$ be the $\frac{37}{2}$-surgery manifold on (-2,3,7)-pretzel knot. Then $\pi_{1}(Y)$ admits no non-cyclic $SO(3)$ representation. Therefore any non-trivial null-homologous knot $K\subset Y$ has non-trivial A-polynomial.
\end{ex}

\subsection{The $S^{3}$ case}

Now let's go back to $S^{3}$ and use the technique in \cite{KM} to prove Theorem \ref{3-sphere case}.
Let $K$ be a non-trivial knot in $S^{3}$. By Remark \ref{reducible}, we can assume that $A_{K}(m,l)=(l-1)f(m,l)$.  By Theorem \ref{non-trivial in sphere}, we have $f(m,l)\neq 1$. Denote the algebraic curve defined by $f(m,l)$ in $\mathds{C}\oplus\mathds{C}$ by $C$.

Consider its intersection with the unit torus $T^{2}=\{(m,l)\mid|m|=|l|=1\}$. It is a closed subset of the torus. The following fact is straightforward by Lemma \ref{dimension} and formula (\ref{boundary}).
\begin{lem}\label{intersection with unit torus}
 If $(\theta,\eta)\in R_{K}$ and $\eta\neq 2k\pi$, then either $(e^{i\theta},e^{i\eta})$ is a $0$-dimensional component of $h^{-1}(\overline{\text{Im}(r)})$ or $(e^{i\theta},e^{i\eta})\in C\cap T^{2}$.
\end{lem}

Now we can prove Theorem \ref{3-sphere case}:

\begin{proof}
Suppose for some unit length $l_{0}=e^{i\eta_{0}}$ ($\eta_{0}\in[0,2\pi)$), we can't find $m_{0}=e^{i\theta_{0}}$ such that $f(m_{0},l_{0})=0$. Since $C\cap T^{2}$ is a compact set, there exists $\epsilon>0$ such that this holds for any $\eta_{0}'\in [\eta_{0}-\epsilon, \eta_{0}+\epsilon]$.  Thus we can assume that $\eta_{0}\neq 0$ and the line $\{l=e^{i\eta_{0}}\}\subset \mathds{C}\oplus\mathds{C}$ does not pass through the $0$-dimensional components of $h^{-1}(\overline{\text{Im}(r)})$. By Lemma \ref{intersection with unit torus}, we see that $R_{K}$ does not intersect the circle $\eta=\eta_{0}$. Thus the translation $R_{K}+(0,-\pi)$ (denote it by $S_{K}$) does not intersect the circle $\eta=\eta_{0}-\pi$. By Definition \ref{pillowcase}, $S_{K}$ is symmetric under the reflection at the origin. Thus $S_{K}$ does not intersect the circle $\eta=\pi-\eta_{0}$. We assume that $\eta_{0}\geq\pi$ for convenience and the other cases can be proved similarly. Now consider the following line segments:
\begin{itemize}
\item $c_{1}=\{(\theta,\eta_{0}-\pi)|\theta\in[0,\pi]\}$;
\item $c_{2}=\{(\theta,\pi-\eta_{0})|\theta\in[-\pi,0]\}$;
\item $c_{3}=\{(0,\eta)|\eta\in [\pi-\eta,\eta_{0}-\pi]$;
\item $c_{4}=\{(\pi,\eta)|\eta\in [0,\eta_{0}-\pi]$
\item $c_{5}=\{(-\pi,\eta)|\eta\in [\pi-\eta_{0},0]\}$
\end{itemize}

 \begin{figure}[htbp]
\centering
\includegraphics[angle=0,width=0.6\textwidth]{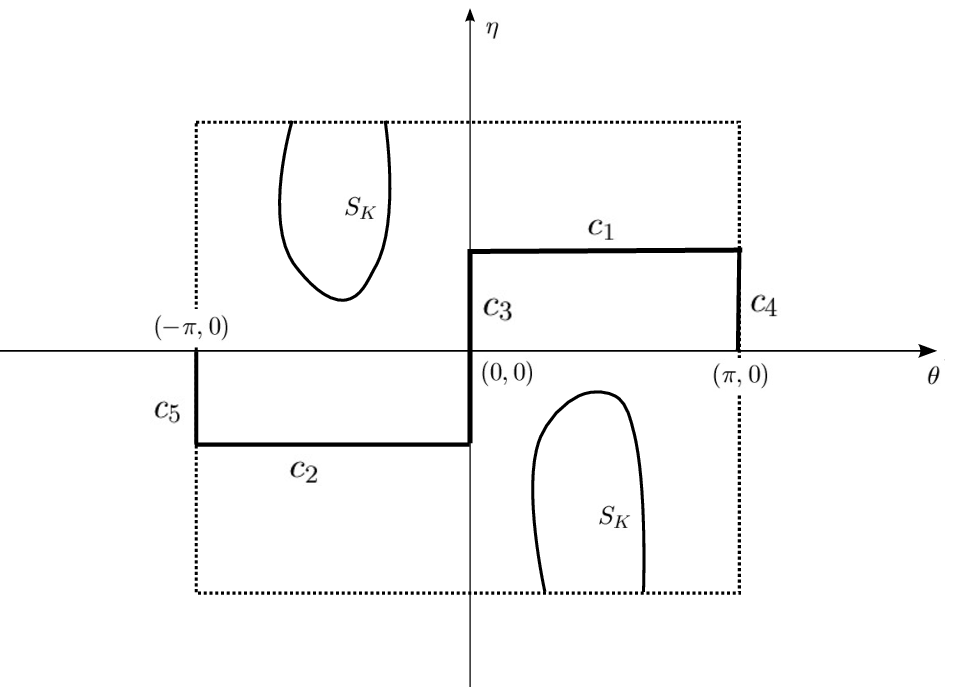}
\caption{} \label{meinv}
\end{figure}

   By the discussion before, $S_{K}$ does not intersect $c_{1}$ and $c_{2}$. Since $S^{3}$ is simply connected, $S_{K}$ does not intersect $c_{3},c_{4},c_{5}$. We can join all these segments to get a piecewise linear path from $(-\pi,0)$ to $(\pi,0)$ and passing through $(0,0)$. By Remark \ref{close subset}, we can find a small neighborhood $N$ of this path such that $S_{K}\cap N=\emptyset$. It is easy to see that we can find an odd, $2\pi$-periodic function $g:[-\pi,\pi]\rightarrow(-\pi,\pi)$ such that the graph of $g$ lies in $N$. Thus $S_{K}$ does not intersect the graph of $g$. We can find an even function $f:[-\pi,\pi]\rightarrow\mathds{R}$ with period $2\pi$ satisfying $f'=g$.

We consider the rank 2 unitary  bundle $E$ over $Y_{0}$ (recall that $Y_{0}$ is the 0-surgery manifold along $K$) with $c_{1}(E)$ the Poincar\'{e} dual of the meridian. We do a holonomy perturbation along the surgery solid torus in $Y_{0}$ using the even function $f$. By Lemma \ref{critical point}, the critical points of the perturbed Chern-Simons functional correspond to the intersection of $S_{K}$ with the graph of $g$, which is empty. (We have $S_{K}$ instead of $R_{K}$ here, because $E$ has non-trivial first Chern class.) So the perturbed Chern-Simons functional has no critical points. This contradicts Corollary \ref{exist critical point}.
\end{proof}

This theorem has the following corollary, which was also proved by Boden in \cite{Boden}:
\begin{cor}
For any non-trivial knot $K\subset S^{3}$, $\text{deg}_{M}A(M,L)\neq 0$.
\end{cor}

\vskip 0.3 truecm

Jianfeng Lin

Department of Mathematics, University of California Los Angeles, Los Angeles, US

juliuslin@math.ucla.edu

\end{document}